\documentclass{amsart}
\usepackage{graphicx, amsmath, amssymb}
\newtheorem{Theorem}{Theorem}[section]

\newtheorem{Lemma}[Theorem]{Lemma}

\theoremstyle{definition}
\newtheorem{Definition}[Theorem]{Definition}
\theoremstyle{remark}
\newtheorem{rem}[Theorem]{Remark}
\numberwithin{equation}{section}

\newcommand{\R}{\mathbb R}
\newcommand{\N}{\mathbb N}
\newcommand{\Z}{\mathbb Z}

\begin{document}

\title{Structure groups and holonomy in infinite dimensions}%
\author{Jean-Pierre Magnot}%
\address{Institut für Angewandte Mathematik \\
Abt. f. Wahrscheinlichkeitstheorie und Mathematische Statistik\\
 Wegelerstr. 6, D-53115 Bonn, GERMANY}%
\email{magnot@wiener.iam.uni-bonn.de}%

\date{12/12/2002}

\begin{abstract}
We give a theorem of reduction of the structure group of a principal bundle $P$ with 
regular structure group $G$. Then, when $G$ is in the classes of regular Lie groups defined by T.Robart in \cite{Rob}, 
we define the closed holonomy group of a connection as the minimal closed Lie subgroup of $G$ for which the previous theorem of reduction can be applied. 
We also prove an infinite dimensional version of the Ambrose-Singer theorem: the Lie algebra of the holonomy group is spanned by the curvature elements. 
\end{abstract}

\maketitle


\section*{Acknowledgements} 
The first investigations on this subject were held during the author's PhD thesis with Professor Sylvie Paycha 
at the University Blaise Pascal, Clermont-Ferrand, France. Most of the results of this paper arose 
at the begining of a post-doctoral fellowship  supported by a Bourse Lavoisier. 
The author would like to thank Professor Sylvie Paycha, 
Professor Sergio Albeverio for inviting him for a postdoctoral fellowship at the Institut fur Angewandte mathematik, Bonn, Germany, 
and also Professor Tilmann Wurzbacher and Professor Claude Roger for inviting him to expose some of these 
results, once at the s\'eminaire de Physique Math\'ematique  of the Institut Girard Desargues, Univ. Lyon I ,France , and also
at the conference ``G\'eom\'etrie en dimension infinie et application \'a la th\'eorie des champs'' held at the CIRM, Marseille, France. 
The author is also indebted to many persons present at these two talks for stimulating questions that influenced the final presentation of this article.

\section*{Introduction}

In the finite dimensional setting, the holonomy group $H$ of a fixed connection $\theta$ with curvature $\Omega$ 
on a principal bundle $P$ with basis $M$ and with structure group $G$ can be defined from two (equivalent) points of view: 

1) $H$ is the smaller Lie group for which there is a principal bundle $Q$ of basis $M$ and with structure group $H$ on which the connection $\theta$ reduces,

2) $H$ is the group of the holonomies of all the horizontal lifts starting at a fixed point in $P $ of a loop based at a fixed point in $M$.

From the first point of view, it enables one to define all Lie groups $H_1$ for which there is a principal bundle $Q_1$ where the connection $\theta$ reduces. 
The second point of view then enables to build $Q_1$ knowing $H_1$. The finite dimensional Ambrose-Singer theorem \cite{as} states that the Lie algebra of the holonomy group is spanned by the curvature elements.
Basic facts about the holonomy in finite dimensions are recalled in section 3. 

In a paper of 1988, Freed \cite{F2} claimed that there was no adequate Ambrose-Singer theorem to deal with 
the Chern form he was defining on loop groups. 
Ten years before,  answering to a question raised by Penot \cite{Pen}, Vassiliou \cite{V} showed 
how far the original 
proof of Ambrose and Singer could be generalized to Banach bundles.  
A crucial ingredient of this proof is that the holonomy group defined from the point of view 2) is already a Lie subgroup of $G$ in the sense of Bourbaki, 
and this paper is also based on a careful application to the Frobenius theorem on Banach vector bundles. 
But this generalization does not fit with Freed's framework. 
Our purpose here is to give an infinite-dimensional generalization of the concept of holonomy group, 
holonomy bundle and of the Ambrose-Singer theorem that fits with Freed's framework.
  
\vskip 12pt

Since we want to consider the framework which is as general as possible, we need to take under consideration 
not only Banach Lie groups, but also some Lie groups $G$ modelled on complete locally convex topological spaces. 
The principal bundle $P$ and the basis manifold $M$ are assumed to be modelled on locally convex topological spaces. 

We remark that there is no adequate Frobenius theorem in this framework. Moreover, to our knowledge, 
there was no proved theorem of reduction for a fixed connection $\theta$. 
Thus, our approach of the problem of defining the Holonomy group and the Holonomy bundle, which is based on 
a definition of the holonomy group with respect to the point of view 1), is the following  : 

- first, we give a theorem to decide whether a connection can reduce while reducing 
the structure group of the principal bundle in section 4. 

- Then, we define \underline{if it exists} the holonomy group $H^{red}$ (resp. the restricted holonomy group $H^{red}_0$) 
as the structure group (resp. the connected component of the structure group) 
of the smaller principal bundle where the connection reduces.

- Finally, using results from Robart \cite{Rob} that we recall in section 1, and remarking that 
we can restrict the existence of the holonomy group to the problem of enlargeability 
of a Lie algebra which is minimal for inclusion, 
we show that we can define  the \textit{closed} restricted 
holonomy group (section 5) and the \textit{closed} holonomy group (section 6) 
when the structure group $G$ belongs to the classes defined in \cite{Rob}. 

The Lie groups $H^{red}$ and $H^{red}_0$ contain the restricted holonomy group and the holonomy group defined by the approach 2), 
but the approaches 1) and 2) do not seem to be equivalent in infinite dimensions. 
However, we do not have any example where the approaches 1) and 2) define different holonomy groups. 
This question remains open, and we make final remarks in chapter 7 that could lead 
to understand better the relationship between these two approaches.

In another paper \cite{Ma}, we will show how these constructions fit to the question about the holonomy in infinite dimensions raised by Freed \cite{F2}.

\section{Enlargeable Lie algebras}

One great problem in the theory of Lie groups and Lie algebras in infinite dimensions 
is to decide whether a Lie algebra is the Lie algebra of a Lie group.
In this article, all the Lie groups are assumed to be modelled on complete locally convex vector spaces. 
We recall in this section some results, that can be found in \cite{Rob}. 
Recall that, in the sense of Bourbaki \cite{Bou}, $H$ is a Lie subgoup of a Lie group $G$ if and only if 
it is a closed subgroup of $G$ and the Lie algebra 
$\frak h$ of $H$ splits the Lie algebra $\frak g$ topologicaly. It is rarely the case. That is why we prefer  the 
following definition of Lie subgroups that is used e.g. in \cite{Rob}:

\begin{Definition}
Let $G$ be a Lie group and $H$ be a subgroup of $G$. Then, $H$ is a Lie subgroup of $G$ if and only if 

- $H$ is a Lie group

- the inclusion $H \subset G$ is a morphism of Lie groups. 
\end{Definition}

We have the analogous definition for Lie subalgebras:

\begin{Definition}
Let $\frak g$ be a Lie algebra and $\frak h$ be an (algebraic) vector subspace of $\frak g$ which is stable under the Lie bracket $[.,.]$ of $\frak g$. 
Then,  $\frak h$ is a Lie subalgebra of $\frak g$ if and only if 

- the topology of $\frak h$ is stronger that the topology of $\frak g$, i.e. the inclusion $\frak h \subset \frak g$ is a morphism of topological vector spaces.

- the Lie bracket $[.,.]$ is continuous as a map $\frak h \times \frak h \rightarrow \frak h$ with respect to the topology of $\frak h$.
  \end{Definition}
Let us first restrict ourselves to the following class of Lie groups. 
\begin{Definition} (\cite{Om}, \cite{Mil}, see e.g. \cite{Rob})

A Lie group $G$ of Lie algebra $\frak g$ is \textit{Omori-Milnor regular} (or \textit{regular} for short) if and only if

\begin{item} 1) for any $v \in C^\infty([0,1],\frak g)$, the differential equation $$g^{-1}{dg(t) \over dt} = v(t)$$ has 
an unique smooth solution $\gamma_v$.
\end{item}
\begin{item} 2) the map $v \mapsto \gamma_v(1)$ is smooth for the topology of uniform comvergence in $ C^\infty([0,1],\frak g)$.
\end{item}
\begin{item} 3) $\gamma_v$ is obtained by product integral in $G$ (see \cite{Om} for the definition of the product integral). \end{item}
\end{Definition}

Actually, we do not know if there exists a Lie group which is not regular. 
On  a regular Lie group, the exponenial map is defined but there are examples 
where it is not locally injective nor locally surjective \cite{Om}, see e.g. \cite{Mil}. 
In this article, we consider only such groups. 
In this context, the tangent space $TG$ is trivial, i.e. $TG \sim G \times \frak g$, 

Let us now give some results about the following problem: let $G$ be a regular Lie group with Lie algebra $\frak g$. 
Let  $ \frak h$ be a Lie subalgebra of $\frak g$. 
When is there a Lie group $H$ of Lie algebra $\frak h$ which is a Lie subgroup of $G$ ? For this, let us give the following definitions:

\begin{Definition}\cite{Rob}

\begin{item} 1) A Lie algebra is called \textit{CBH} if and only if its Campbell-Baker-Hausdorff serie 
$$ (v,w) \in U \times U \mapsto h_v(w)=v + w + {1 \over 2}[v,w] + {1 \over 12}\Big( [v,[v,w]]-[w,[v,w]]\Big) +...$$ 
is analytic on $U \times U$, where $U$ is a neighborhood of $0$ in $\frak g$, that is if and only if the Campbell-Baker-Hausdorff serie converges for any $(v,w) \in U \times U$.

A Lie group is called \textit{CBH} if and only if its Lie algebra is CBH.
\end{item}
\begin{item}
2 ) A Lie group $G$ is of the \textit{first type} if and only if the exponential map furnishes a chart at the neighborhood of the identity.
\end{item}
\begin{item}
3 ) A Lie group $G$ of Lie algebra $\frak g$ is of the \textit{second type} if and only if  $\frak g$ decomposes as a sum of Lie algebras $$ \frak g = \oplus_{i=1}^m \frak g_m$$ and the exponential map furnishes a chart at the neighborhood of the identity by the map $$ (v_1,...v_m) \in \frak g_1 \times ... \times \frak g_m \mapsto \prod_{i=1}^m exp(v_i) \in G.$$
\end{item}
\end{Definition}

A CBH Lie group is of the first type, and a group of the first type is of the second type. 
For examples of CBH Lie groups, see \cite{Rob} and also \cite{Gl}.

\vskip 10pt
Let us now give the results of \cite{Rob} we shall use in this article. 
Recall that, if $\frak h$ is a Lie subalgebra of $\frak g$, we do not assume a priori that $\frak h$ is closed in $\frak g$, 
nor the existence of a topological splitting. 

\begin{Theorem}\label{LieII}
Let $\frak h$ be a Lie subalgebra of $\frak g$, and assume that $\frak g$ is the Lie algebra of a Lie group $G$. In the following cases:
\begin{item} 1) If $G$ is CBH and $\frak h $ is closed in $\frak g$,
\end{item}
\begin{item} 2) If $G$ is CBH and $\frak h $ is CBH, 
\end{item}
\begin{item} 3) If $G$ is of the first type and $\frak h $ is closed in $\frak g$,
\end{item}
\begin{item} 4) If $G$ is of the second type and $\frak h $ is closed in $\frak g$,
\end{item}
\vskip 5pt
\noindent
there is an unique connected Lie subgroup $H$ of $G$, of the same type as $G$, that has $\frak h$ as Lie algebra.
\end{Theorem}

Let us now give the following (easy but useful) lemma:
\begin{Lemma}
Let $\frak g$ be a Lie algebra. Let $X\subset \frak g$ be a topological space, with continuous inclusion. Let $A$ be a set of the closed subalgebras $\frak h \subset \frak g$ such that 
$$ X \subset \frak h \subset \frak g.$$ 
Then, $ A$ has an unique minimal element for the inclusion.
\end{Lemma}
The proof of this lemma is straightforward, taking as minimal element $$\frak h_0 = \cap_{\frak h \in A} \frak h,$$
which is a closed subalgebra of $\frak g$.

\section{Structure groups}

In all the article, $P$ is a principal bundle, with basis $M$ and with structure Lie group $G$, 
of Lie algebra $\frak g = T_eG.$ Let $\pi$ be the fiberwise projection from $P$ on $M$. We do not specify 
first the model spaces of these manifolds. They can be Hilbert, Banach, Frechet, locally convex topological spaces, 
or even take place in the "convenient setting" \cite{KM}. 
All the differentiable structures are assumed to be smooth. 
We recall that a connection $\theta$ on the principal bundle $P$ is a vector bundle projection on the vertical bundle
$TP \rightarrow VTP$, and can also be seen as a 1-form on $TP$ with values in $G$ which is $Ad$-invariant. Its curvature, $\Omega$, 
can be also viewed from these two points of view. Recall that if we have an Ad-invariant  map $f : P \rightarrow \frak g$ ( i.e. such that for any $p \in P$ and for any $g \in G$, $f(p.g) = Ad_{g^{-1}}f(p)$), 
then the induced vertical vector field $p \mapsto f(p)^*$ of  $TP$ is invariant under the action of $G$ on $TP$. 
For basic facts on principal bundles, we refer to \cite{KM}.

\vskip 10pt
1) On a principal bundle $P$, using a Lie group morphism $\rho: G \rightarrow G'$,
we build the principal bundle $P\times_\rho G'$. Then, the connection $\theta$ defines a connection $\theta'$ 
on $P\times_\rho G'$, see e.g. \cite{KMS}, p. 107-108, or \cite{KM} p. 383. 

2) The \textit{reduction:} Let $G_1$ be a Lie group with Lie algebra $\frak g _1$, and $\rho : G_1   \rightarrow G$ 
be an injective morphism of Lie groups. When we write $\rho(G_1)$ instead of $G_1$, 
we consider the group endowed with the topology and the differentiable structure induced by $G$, 
which can be weaker than the initial topology and differentiable structure $G_1$. 
 
If there exists a family of local trivializations of $P$ whose transition functions take their values in $\rho(G_1)$, and that are smooth with respect to the topology of $G_1$ 
then we can build a $G_1$-bundle $P_1$ 
with respect to these trivializations. 
$P_1$ is called a \textit{ reduction}  of $P$ and there is an inclusion $i : P_1 \rightarrow P,$ 
which is a morphism of principal bundles.  
The connection $\theta$ \textit{reduces} to $P_1$ if and only if there is a connection 
$\theta_1$ om $P_1$ that is the pull-back of the connection $\theta$ on $P$. 
This is the case when, on each local trivialization of $P$ defined before, 
the connection $\theta$ reads as a smooth local 1-form on $M$ with values in $\frak g _1$.

\vskip 10pt

As one can easily see, the first way of changing the structure group of $P$ is possible in any case. 
The second one, the reduction, requires more attention, and one of the aims of this article is to provide a tool 
to decide whether reduction is possible in an infinite dimensional setting, 
and then to extend the finite dimensional Ambrose-Singer theorem \cite{as}. 

In order to prove the Ambrose-Singer theorem, we shall also need to change the basis manifold, 
which is possible by \textit{pull back}, see e.g. \cite{KM}, p. 377.: 
Let $N$ be a manifold and let $f:N \rightarrow M$ be a smooth map. 
Then, one defines $f^*P$, the pull back of the principal bundle $P$ by $f$, 
which is a principal bundle with structure group $G$ and of basis $N$. Then, a connection $\theta$ always defines a connection $f^*\theta$ on $f^*P$. We use this construction when $f$ is the projection from the universal covering $N= \Sigma M$ of $M$ to $M$. In that case, we write $\tilde \theta$ instead of $f^
* \theta$, and the curvature of $\tilde \theta$ is the same as the curvature of $\theta$.
\section{Holonomy groups}
Let us \textbf{fix} now a connection $\theta$ on the principal bundle $
P$, with regular structure group $G$ modelled 
on a complete locally convex topological vector space and of basis $M$. 
We assume that 
there exists in $M$ a dense countable subset $Q$, that $M$ is modelled on a locally convex vector space, 
 and  that the fundamental group $\pi_1(M)$ is countable. 

\subsection{Holonomy induced by horizontal paths}
Let $p \in P$ be a fixed point, and let $x = \pi(p)$ be the basepoint of $p$. 
Let $\gamma$ be a piecewise smooth path on $M$ starting at $x$, i.e. 
$\gamma: [0,1] \rightarrow M$. Then, \cite{KM}, the differential equation on $p$
\begin{eqnarray*}  c(0) & = & p \\
\pi(c(t)) & = & \gamma(t) \hbox{ for } t \in [0,1] \\
\theta(dc(t) / dt) & = & 0  \end{eqnarray*}
has an unique solution called the horizontal lift of $\gamma$ in $P$, 
that we note by $H\gamma$. 
We define the Holonomy group from the point of view 2) of the introduction:
$$H^{curv}(p) = \{ g \in G \, | \, H\gamma(1)=p.g, \,  \gamma(1) = x \}.$$
Given two points $p,p' \in P$, the holonomy groups $H^{curv}(p) $ and $H^{curv}(p')$ are conjugate in $G$. If $p$ and $p'$ can be joint by a horizontal path, $H^{curv}(p) = H^{curv}(p')$, 
and if $P$ is finite dimensional, 
$H^{curv}(p)$ is a Lie group, which appears to be the smaller subgroup with which 
one can reduce the structure group of $P$ and also reduce the connection $\theta$. 
Note that they are often not closed Lie subgroups of $G$, as shows this example:
\vskip 10pt
\noindent
\textbf{Example:} Let $P=S^1\times S^1=\R^2 / \Z^2$ be the trivial principal bundle 
over $S^1$ with structure group $S^1$. We identify each tangent space of 
$P$ with $\R^2$ in the canonical way, the second coordinate being the vertical coordinate. 
Let us now consider the connection $\theta$ defined by $\theta : (x,y) \in \R^2 \mapsto \theta(x,y) = y + x.\pi$.
 Then, the corresponding holonomy group is given by $\Z[\pi] / \Z$, which is a countable dense subgroup of 
$S^1 = \R/\Z$. The holonomy group is itself a Lie group (of dimension $0$), but not a closed Lie subgroup of $S^1$.
\vskip 10pt
In order to avoid these problems, one can also consider the restricted holonomy group, 
the connected component of the neutral element of the holonomy group, 
which is obtained considering only the loops in $x$ that are homotopy trivial. 
From another point of view, the restricted holonomy group is also the holonomy group of the connection $\tilde \theta$ induced by $\theta$
on the pull back of $P$ on the universal covering $\Sigma M$ of $M$. 
This pull-back can also be done in infinite dimensions.

\subsection{Holonomy, curvature and Ambrose-Singer theorems: known
results}

We define the horizontal bundle as $HTP = Ker(\theta)$. This bundle gives a decomposition 
$TP = HTP \oplus VTP$. The curvature $\Omega$ of $\theta$ measures the vertical component of the 
bracket of two sections $X,Y$ of the horizontal bundle, 
that is, given $p \in P$, $\Omega_p(X,Y) = \theta_p([X,Y])$. 
One has to be careful because the Lie algebra spanned by the curvature elements 
needs not to have the same dimension at each point $p \in P$ (see e,g. \cite{Li}, p. 141). 

Nevertheless, we have the finite dimensional Ambrose-Singer theorem:

\begin{Theorem} 
\cite{as}
The Lie algebra of the restricted holonomy group 
at $p \in P$ is spanned by the curvature elements.
\end{Theorem}

Let us precise the following: by ``curvature elements'', we understand the curvature elements $\Omega_{p'}(X,Y)$, where $p' \in P$ can be joint to $p$ by a horizontal path, and $X,Y \in T_{\pi(p')}M$. 
The proof of this theorem uses two tools : the first one is 
that the restricted holonomy group is finite dimensional. 
The second one is the Frobenius theorem for vector bundles. 
Such a theorem was not yet proved in infinite dimension, 
except on a particular case of holonomy group embedded in a Banach Lie group in \cite{V}. 
In this article, the very strong assumptions on the holonomy group replaces the finite dimensional setting, 
which leads the author to apply very easily the Frobenius theorem for Banach vector bundles. 

\vskip 10pt

In our setting, there is no Frobenius theorem, and we want to consider all the connections that exist on $P$. 
This is the reason why our approach is quite different.  

\section{Reduction of the structure group}
Let us first consider the following question: given a Lie subgroup $G_1$, with Lie algebra $\frak g_1$, of $G$, 
when is it possible to reduce the bundle $P$ to a principal bundle $P_1$ with structure group $G_1$ ? 
The most famous condition, in finite dimension and if the basis manifold is simply connected, 
is the existence of a connection with $\frak g_1$-valued curvature. 
This result is of course linked with the finite-dimensional proof of the Ambrose-Singer theorem. 
We prove here an infinite dimensional version of the theorem of reduction. 
This allows us, in the next chapter, to deal with an infinite dimensional construction of the holonomy groups from the point of view 1) of the introduction. In the following, $p$ is a fixed point in $P$ with basepoint $x$ in $M$.

\begin{Theorem} \label{Courbure}
We assume that $G_1$ and $G$ are regular Lie groups and that $G_1$ is modelled on a complete locally convex topological vector space. 
Let $\rho: G_1 \mapsto G$ be an injective morphism of Lie groups, 
if there exists $\theta$ a connection on $P$, with curvature $\Omega$, such that, for any smooth 1-parameter family $Hc_t$ of horizontal paths starting at $p$, for any smooth vector fields $X,Y$ in $TM$,  
\begin{eqnarray} s, t \in [0,1]^2 & \rightarrow & \Omega_{c_t(s)}(X,Y)  \label{g1}\end{eqnarray} 
is a smooth $\frak g_1$-valued map,
\noindent
and if $M$ is simply connected, then it is possible to reduce the structure group $G$ of $P$ to $G_1$.
\end{Theorem}  

We shall of course try to apply this theorem to find an infinite dimensional version of the Ambrose-Singer theorem. Let us first give a remark that makes easier the application of this theorem in concrete situations: 

\begin{rem} If $\frak g_1$ is stable under the adjoint action of $G$, the map \ref{g1} is a smooth $\frak g_1$-valued map for any 1-parameter family $Hc_t$ of horizontal paths if and only if on any local trivialization of $P$, the curvature $\Omega$ reads as a smooth $\frak g_1$-valued 2-form. Otherwise, it is sufficient to find a family of local trivializations of $P$ induced by an open covering of $M$ for which the curvature $\Omega$ reads as a smooth $\frak g_1$-valued 2-form.    \end{rem}
 
Before giving the proof of this theorem, we need three lemmas, 
which are well-known results in finite dimensions (see e.g. \cite{Li}). Following \cite{KM}, 
if $ C^{\infty}_x([0,1],M) $ is the set of smooth paths on $M$ starting at $x$,
we note $$Pt : C^{\infty}_x([0,1],M) \times [0,1] \times \pi^{-1}(x) \rightarrow P$$ the parallel transport with respect to $\theta$, which is a smooth map. 

\begin{Lemma} \label{cartes}
Let $p \in P$. Let $x$ be the basepoint of $p$. 
Let $U$ be a star-shaped neighborhood of $x$, 
that we idendify in sake of simplicity for this lemma with a  
star-shaped neighborhood of 0 in the model space of $M$. 
We assume that $U$ is contained in the domain $V$ of a smooth trivialization 
$\varphi : P_{|V} \rightarrow V \times G$ of $P$ 
Let $u \in U$ and $t \in [0,1]$. We define $f(u,t)= tu \in U$. 

Let \begin{eqnarray*}
\psi : U & \rightarrow & P_{|U} \\
     u & \mapsto & Pt(f(u,.), 1  , p) \end{eqnarray*}

Let \begin{eqnarray*} \Psi : U \times G & \rightarrow & P_{|U}\\
(u,g) & \mapsto & \psi(u).g \end{eqnarray*}

and $$ \tilde\Psi = \Psi \circ (Id_{U} \times \rho).$$

Then, $\Psi$ is a local trivialization of $P$ over $U$. Moreover, $\theta \circ D\psi$ is a smooth 
$\frak g_1$-valued form on $U$. 
\end{Lemma}

We need now to know how the horizontal lifts of paths behave in this trivialization.

\begin{Lemma}\label{horizontal}
We use here the notations of the two last lemmas.
We assume also that $U$ is convex.
\begin{item}
(i) Given a path $\alpha: [0,1] \rightarrow U$ starting at $x$, if $H\alpha$ is its horizintal lift starting at $p$, we have $H\alpha(1) \in \Psi ( U \times \rho(G_1))$, and there exists a path $H\alpha_1 :[0,1] \rightarrow U \times G_1$ such that $H\alpha=\tilde\Psi  \circ H\alpha_1.$
\end{item}
\begin{item}
(ii) Let $h:[0,1]^2 \rightarrow U $ be an homotopy equivalence between two paths $h(0,.)$ and $h(1,.)$ starting at $x$ and finishing in $U$. Let $Hh(0,.)$ and $Hh(1,.)$ be their horizontal lifts starting at $p$. Then, there is $g_1 \in G_1$ such that $\tilde\Psi^{-1}(Hh(0,1)) = \tilde\Psi^{-1}(Hh(1,1)).\rho(g_1).$    
\end{item}
\end{Lemma}

Then, the following lemma will be useful when dealing with homotopy:

\begin{Lemma}\label{final}
Let $\alpha$ and $\beta$ be two paths on $U$. Let $q_\alpha \in \pi^{-1}(\alpha(0))$ and $q_\beta \in \pi^{-1}(\beta(0))$. Let $H\alpha$ and $H\beta$ be the horizotal lifts of $\alpha$ and $\beta$ starting at  $q_\alpha$ and $q_\beta$. We set $\tilde \Psi^{-1}\circ H\alpha = (\alpha,\gamma_\alpha)$ and $\Psi^{-1}\circ H\beta = (\beta,\gamma_\beta)$. 

Let $g = \gamma_\beta^{-1}(0).\gamma_\alpha(0)$, with $^{-1}$ as 
the inverse map in $G$. Then, for any $t \in [0,1]$, there exists $g_1(t), g'_1 \in G_1$ such that 
$\gamma_\beta(t) = \gamma_\alpha(t).g_1^{-1}(t).g.g_1'(t)$. Moreover, 
the maps $t \mapsto g_1(t)$ and $t \mapsto g'_1(t)$ are smooth in $G_1$.
\end{Lemma}
\noindent 
Let us now give the proofs of the three lemmas, and then the proof of Theorem \ref{Courbure}:
\vskip 10pt
\noindent
\textbf{Proof of Lemma \ref{cartes}:}

We already know that $\Psi$ is a smooth map, since $Pt$ is smooth. We also know that $\psi$ is
 a smooth injective map, and that $\pi \circ \psi = Id_{U}$. 
Hence, $\Psi$ is a smooth local trivialization of $P$ over $U$. 

Let us calculate $D\Psi^{-1} \circ \theta \circ D\psi$.
Let $c:]-\epsilon, \epsilon[ \rightarrow U$ be a smooth path such that $c(0) = u \in U$. Let $h(t,s) = f(c(s),t)$. 
Let $\tilde \theta$ be the pull-back of $\theta$ by $\Psi$ on $U$.

Then, following the proof of the claim of \cite{KM}, p.424, with a connection with non vanishing curvature, 
we have:

\begin{eqnarray*} \partial_s(h^*\tilde \theta)(\partial_t) & = &  \partial_t(h^*\tilde \theta)(\partial_s) - d(h^*\tilde \theta)(\partial_t,\partial_s) - (h^*\tilde \theta)([\partial_t,\partial_s]) \\
& = & \partial_t(h^*\tilde \theta)(\partial_s) - d(h^*\tilde \theta)(\partial_t,\partial_s)\\
& = & \partial_t(h^*\tilde \theta)(\partial_s) + ad_{(h^*\tilde \theta)(\partial_t)}((h^*\tilde \theta)(\partial_s)) - (\Psi^*\Omega)(h_*(\partial_s),h_*(\partial_t)) \\
& = & \partial_t(h^*\tilde \theta)(\partial_s) - (\Psi^*\Omega)(h_*(\partial_s),h_*(\partial_t)) \quad \hbox{since } (h^*\tilde \theta)(\partial_t)= 0 .\end{eqnarray*}

Then, 

$$\partial_s(\Psi \circ \psi \circ c) = (\partial_s c(s), \partial_s \tilde\gamma(1,s)),$$
remarking that $(u,e) = \Psi^{-1} \circ \psi(u)$. 
We now calculate $\partial_s \tilde\gamma(1,s)$, 
\begin{eqnarray*}
\partial_s \tilde\gamma(1,s) & = & \int_0^1 \left( \partial_s(h^*\tilde \theta)(\partial_t) \right)(t
) dt \\
& = & \int_0^1 \Big(\partial_t(h^*\tilde \theta)(\partial_s) 
- (\Psi^*\Omega)(h_*(\partial_s),h_*(\partial_t))\Big)(t) dt
\\
& = & (h^*\tilde \theta)(\partial_s)(1,s) - \int_0^1 (\Psi^*\Omega)(h_*(\partial_s),h_*(\partial_t))(t) dt\end{eqnarray*}
Finally, we have: 
\begin{eqnarray} D\Psi^{-1} \circ \theta(\partial_s(\psi \circ c)) &=& (h^* \tilde \theta)(\partial_s h(1,s),\partial_s \gamma(1,s))\\
& = & \int_0^1 \Big(\Omega(h_*(\partial_s),h_*(\partial_t))(t) \Big)dt \label{evolution}
 \end{eqnarray}
Since $\frak g_1$ is complete, this integral exists and belongs to $\frak g_1$.
\qed
\vskip 10pt
\noindent
\textbf{Proof of the Lemma \ref{horizontal}:} 

(i) We have that $\theta(\partial_s(\psi \circ \alpha))$ is an integral on the curvature elements (see the proof of the last lemma). 
Looking at this result more precisely, reparametrizing equation \ref{evolution}, setting $c=\alpha$, we have that
\begin{eqnarray}\theta(\partial_s(h(s,t)) & = & \int_0^t \Big(\Omega(h_*(\partial_s),h_*(\partial_t))(u) \Big)du,\end{eqnarray} and hence that $$\partial_t\Big(\theta(\partial_s(h(s,t))\Big) = \Omega(h_*(\partial_s),h_*(\partial_t).$$
Recall that $ \rho^{-1}\circ\Omega(h_*(\partial_s),h_*(\partial_t)) $is smooth.
Then, integrating this equality in $G_1$ instead of $G$, we get a path $\alpha_1$ in $U \times G_1$. Then we consider the following differential equation, that defines $H\alpha_1$: 
$$ \left\{ \begin{array}{ccl}H\alpha_1(0)&=&e\\
 H\alpha_1(t) & = & (\alpha(t),\gamma(t)) \in U \times G_1\\
D\rho\big( \gamma(t)^{-1}\partial_t\gamma(t) \big) &=& Ad_{(\rho \circ \gamma^{-1})(t)} \Big( \theta \circ D\Psi\circ(Id \times\rho)(\partial_t\alpha_1)(t)\Big) , \end{array} \right.$$
setting  $H\alpha=\Psi \circ (Id_U \times \rho) \circ H\alpha_1,$ we get (i).

(ii) comes easily from the continuity of the horizontal lift of paths, using the fact that $U$ is contractible, and applying (i) to the path 
\begin{equation} \label{hg1} c(t)= \left\{\begin{array}{ccl} Hh(0,3t) & \hbox{if} & t \in [0,1/3] \\
   Hh(3t-1,1)& \hbox{if} & t \in [1/3,2/3]\\
 Hh(1,3-3t)& \hbox{if} & t \in [2/3,1].\end{array}\right.\end{equation}
\qed
\vskip 10pt
\noindent
\textbf{Proof of Lemma \ref{final}:} Reparametrizing Lemma \ref{horizontal}, (i), we have that, for $t \in [0,1]$, there exists $g_1(t), g'_1(t) \in G_1$ such that $\gamma_\alpha(t) = \gamma_\alpha(0).g_1(t)$ and $\gamma_\beta(t) = \gamma_\beta(0).g'_1(t)$. Then, we get $\gamma_\beta(t) = \gamma_\alpha(t).g_1^{-1}(t).g.g_1'(t)$. By equation \ref{hg1}, we have that the paths $t \mapsto g_1(t)$ and $t \mapsto g'_1(t)$ are smooth paths with values in $G_1$. \qed   
\vskip 10pt
\noindent
\textbf{Proof of Theorem \ref{Courbure}:}
We consider the covering $\{U_x\}$ of $M$ indexed by the points $x \in M$ such that each $U_x$ is star-shaped. Recall that we have assumed that $M$ has a countable dense subset $Q$.

From the map $\psi_{p}$, we define as in Lemma \ref{horizontal} the map $\Psi_{p} : U_x \times G \rightarrow P$ such that $\Psi_{p}(u,g)=\psi_{p}(u).g$ and $\tilde\Psi_{p} : U_x \times G_1 \rightarrow P$ such that $\tilde\Psi_{p}(u,g_1) = \Psi_{p}(u,\rho(g_1))$. We define $p_0 = p$

Let us now build by induction a sequence of countable families of local trivilizations of $P$.  

\textit{Rank 0:} We have built the map $\Psi_0=\Psi_p$, $\tilde\Psi_0=\tilde\Psi_p$ and set $U_0=U_x$. We define $p_0 = p$.
 
\textit{Assume that we have built the maps at the rank $n-1$ for some $n \in \N^*$. Let us build the family of rank $n$:} The maps $\Psi_{a}$, $\tilde\Psi_{a}$, the points $p_a$ and the sets $U_{a}$ are indexed by the multi-indexes $a \in \N^{n}$ such that $a_1=0$. Let us fix such a multi-index $a$. $Q\cap U_a$ is countable. Then, we order $Q\cap U_a$ to get a sequence $\{q_k\}_{k \in \N}$. Fixing $k \in \N$, we define $b\in \N^{n+1}$, $p_b$, $U_b$ and $\Psi_b$ the following way:

\noindent
- $ b =  (a_1,...,a_{n},k) $

\noindent
- $ p_b = \Psi_a(q_k,e) $

\noindent
- $U_b = U_{q_k} $

\noindent
- $\Psi_b  \hbox{is the map } \Psi \hbox{ of Lemma \ref{horizontal} centered at } p = q_k $

\noindent
- $\tilde \Psi _b = \Psi_b \circ(Id_{U_b} \times \rho)$.

We have built by induction a family of local trivializations indexed on finite sequences of $\N$. Now, we have to check that we get the desired family of local trivializations of $P$, namely, that the family $\{\tilde \Psi_b \}$ has transition maps with values in $G_1$. 

\vskip 10pt

Let $a$ and $b$ be two finite sequences of $\N$ of order $m$ and $n$ such that $U_{a} \cap U_b \neq \not 0$, 
(recall that the transition function $\tau_{a,b}$ is defined by
\begin{eqnarray*} \Psi_{b}^{- 1} \circ \Psi_{a} : U_{a} \cap U_b \times G & \rightarrow & U_{a} \cap U_b \times G  \\
(u,g) &\mapsto &(u,\tau_{a,b}.g)\quad ).\end{eqnarray*}
 We must check that $ \rho^{-1} \circ \tau_{a,b}$ smooth as a $G_1$-valued function. We assume that $m \leq n$. 
Let $z$ be a fixed point of $U_{a} \cap U_{b}$. 
Let us consider the following two paths, starting at $x=\pi(p_0)$ and finishing at $z$. In order to build then, we need also to define the charts $\xi_c : U_c \rightarrow F$  which identify $U_c$  with an open convex neighborhood of $0$ of $F$ for any finite sequence $c$ of $\N$ (this notation was hidden in the three previous lemmas because we were working on a fixed local chaart where as we need here to precise the chart used). We set $$
 \begin{array}{ccl}
\beta_a(t) & = & \xi_0^{-1}(mt\xi_0(p_{(0,a_2)})\hbox{ for } t\in[0,1/m] \\
\beta_a(t) & = & \xi_{(0,a_2)}^{-1}((mt-1)\xi_{(0,a_2)}(p_{(0,a_2,a_3)}) \hbox{ for } t\in[1/m,2/m]\\
& (...) & \\
\beta_a(t) & = & \xi_{(0,a_2,...a_{m-1})}^{-1}((mt-m+2	) \xi_{(0,a_1,...a_{m-1})}(p_{(0,a_1,...a_{m})}) \hbox{ for } t\in[(m-2)/m,(m-1)/m] \\
\beta_a(t) & = & \xi_{(0,a_1,...a_{m})}^{-1}((mt-m+1	) \xi_{(0,a_1,...a_{m})}(z) \hbox{ for } t\in[(m-1)/m,1] \\
\end{array} $$
and we build the same way $\beta_b$, replacing $m$ by $n$ and $a$ by $b$. We define  also the paths :
$$
 \begin{array}{ccl}
\alpha_a(t) & = & \Psi_0(\beta_a(t),e) \hbox{ for } t\in[0,1/m] \\
\alpha_a(t) & = & \Psi_{(0,a_2)}(\beta_a(t),e) \hbox{ for } t\in[1/m,2/m]\\
& (...) & \\
\alpha_a(t) & = & \Psi_{(0,a_2,...a_{m-1})}(\beta_a(t),e) \hbox{ for } t\in[(m-2)/m,(m-1)/m] \\
\alpha_a(t) & = & \Psi_{(0,a_1,...a_{m})}(\beta_a(t),e) \hbox{ for } t\in[(m-1)/m,1] \\
\end{array} $$ 
and $\alpha_b$ replacing $m$ by $n$ and $a$ by $b$. These are piecewise smooth paths with values in $p$, ending at $p_a$ and $p_b$. 
Recall that $M$ is assumed to be connected and simply connected. Hence, there is an homotopy among piecewise smooth paths $h : [0,1]^2 \rightarrow M$ such that $h(0,.)=\beta_a$, $h(1,.)=\beta_b$, $h(.,0)=\pi(p)$, $h(.,1)=z$.
{The horizontal lifts of $\beta_a$ and $\beta_b$ starting at $p$ are the paths $\alpha_a$ and $\alpha_b$}. Since $[0,1]^2$ is compact, there exists $M,N \in \N$ such that, for each $i\in \{0,...,M-1\}$ and $j \in \{0,...,N-1\}$, there is a multiindex $c(i,j)$ (of any order) such that $h([i/M,(i+1)/M]\times [j/N,(j+1)/N])\subset U_{c(i,j)}$. 

Then, we apply the lemma \ref{final} the paths $t \in [0,1] \mapsto Hh(i/M,(j/N)+(t/N))$ and $t \in [0,1] \mapsto Hh((i+1)/M,(j/N)+(t/N))$ in $U_{c(i,j)}$ for each $i\in \{0,...,M-1\}$ and $j \in \{0,...,N-1\}$the following way:

\vskip 5pt
- We first apply it to $t \mapsto Hh(0/M,(0/N)+(t/N))$ and $t \mapsto Hh(1/M,(0/N)+(t/N))$, and then successively to each couple of path $t \mapsto Hh(0/M,(j/N)+(t/N))$ and $t \mapsto Hh(1/M,(j/N)+(t/N))$ for $j = 1,...,N-1$ to have that there exists $g_1^{(0)} \in G_1$ such that $Hh(1/M,1)=p_a.g_1^{(0)}$. We have also proved that, for any $j \in \{1,...,N-1\}$, there exists $g_1^{(0,j)} \in G_1$ such that $Hh(1/M,(j/N)) = Hh(0/M,(j/N)).g_1^{(0,j)}$

\vskip 5pt
- We assume that, for fixed $i$,  there exists $g_1^{(i-1)} \in G_1$ such that $Hh(i/M,1)=p_a.\rho(g_1^{(i-1)})$ and that for any $j \in \{1,...,N-1\}$, there exists $g_1^{(i-1,j)} \in G_1$ such that $Hh(i/M,(j/N)) = Hh(0/M,(j/N)).g_1^{(i-1,j)}$. 

Then, applying the same procedure as before, we consider each couple of path 
 $t \mapsto Hh(i/M,(j/N)+(t/N))$ and $t \mapsto Hh((i+1)/M,(j/N)+(t/N))$ for $j = 0,...,N-1$, on which we apply Lemma \ref{final}. Then, there exists $g_1^{(0,j)} \in G_1$ such that $Hh((i+1)/M,(j/N)) = Hh(0/M,(j/N)).g_1^{(i,j)}$ and $g_1^{(i)} \in G_1$ such that $Hh((i+1)/M,1)=p_a.g_1^{(i)}$. 

\vskip 5pt

We have proved by induction that $$p_b = Hh(1,1) = H(0,1).g_1^{(M-1,N-1)} = p_a.g_1^{(M-1,N-1)}.$$
From Lemma \ref{horizontal}, (i), we have also that the map $z \mapsto \rho^{-1} \circ \tau_{a,b}$ is smooth, remarking that, for $s < (N-1)/N$, $Hh(0,s)$ and $Hh(1,s)$ do not depend on the choice of $z$. \qed

\vskip 10pt
Let us now deal with the case where $M$ is not simply connected.

\begin{Theorem} \label{CourbureII}   
Let $G_1$ a Lie group with Lie algebra $\frak g_1$, and a free representation $\rho_1: G_1 \rightarrow G$. We assume that

\noindent
- for any smooth 1-parameter family $Hc_t$ of horizontal paths starting at $p$, for any smooth vector fields $X,Y$ in $TM$,  
\begin{equation*} s, t \in [0,1]^2 \rightarrow \Omega_{c_t(s)}(X,Y)  \end{equation*} 
is a smooth $\frak g_1$-valued map,

\noindent - for any  horizontal path $H\alpha$ induced by a loop $\alpha$ in $P$, there is $g_1 \in G_1$ such that $H\alpha(1) = H\alpha (0).\rho_1(g_1)$.
\vskip 5pt \noindent
Then, we can reduce the structure group of $P$ to $G_1$.
 
\end{Theorem}

\textbf{Proof of the Theorem \ref{CourbureII}:}
The proof begins like the one of theorem \ref{Courbure}, until we build the family of local trivializations $\Psi_a$, where $a$ is a finite sequence of elements of $\N$ such that $a_1=0$. But here, the additional condition that we have, namely, the holonomy of each horizontal path lies in $\rho_1(G_1)$, gives us directely the result. \qed

\vskip 10pt
Note that, in these two theorems, the connection $\theta$ reduces to the reduced bundle we have obtained. Let us now give the following lemma, which will be useful in what follows, if $M$ is not simply connected.

\begin{Lemma} \label{semidirect}

Let $G_1$ be a connected regular Lie subgroup of $G$ with Lie algebra $\frak g_1$  such that there is a connection $\theta$ with $\frak g_1$-valued curvature.  Assume that its Lie algebra $\frak g_1$ is stable under the adjoint action of $H^{curv}$. Then, there exists a regular Lie group $G'_1$ with connected component $G_1$ for which the theorem of reduction \ref{CourbureII} can be applied.

\end{Lemma}

\noindent
\textbf{Proof of Lemma \ref{semidirect}:}
As we already mentioned, if $\Sigma M$ is the universal covering of $M$ and $\tilde P$ (resp. $\tilde \theta$) is the pull back of $P$ (resp. $\theta$) over $\Sigma M$, $H^{curv}_0(\theta) = H^{curv}(\tilde \theta)$. Since $\Sigma M$ is simply connected, Theorem \ref{Courbure} applies to $\tilde \theta$, and as a consequence $H^{curv}_0(\theta) = H^{curv}(\tilde \theta) \subset G_1$ with continuous inclusion. Thus, in order to show the Lemma, we show that it is possible to make a semi-direct product of $G_1$ with $H^{curv}(\theta) / H^{curv}_0(\theta)$ with is a Lie subgroup of $G$. Recall that, since 
$H^{curv}/H^{curv}_0$  is a group and we have a continuous surjection $\pi_1(M) \rightarrow H^{curv}/H^{curv}_0$, and hence $H^{curv}/H^{curv}_0$ is discrete. 
Let $\{c_i\}$ be a family of loops generating $\pi_1(M)$, and let $g_i$ be their holonomies. Let $g \in G_1$, then, there exists a smooth path  $t \in [0,1]\mapsto g(t) \in G_1$ starting at the neutral element such that $g = g(1)$. Let $v(t) = g^{-1}(t)\partial_tg(t) \in C^\infty([0,1],\frak g_1)$ its logarithmic derivative. $g_i.g.g_i^{-1}$ is the endpoint of the path $t \mapsto g_i.g (t) .g_i^{-1}$, with logarithmic derivative $t \mapsto Ad_{g_i}.v (t)$. 

\begin{equation*} \left.\begin{array}{cl} \forall t \in [0,1], & Ad_{g_i}.v (t) \in \frak g_1 \\
G_1 & \hbox{is regular} \end{array} \right\} \Rightarrow g_i.g.g_i^{-1} \in G_1 \quad \qed \end{equation*}

\section{The restricted holonomy group}

The holonomy group $H^{curv}_0$ induced by horizontal paths is not actually an infinite dimensional Lie group, principaly because its structure is induced by the structure of the non associative monoid of the loops based at $x$.
Howver, using the theorem of reduction of the previous chapter, the definition of the resricted holonomy group as a Lie group turns out to be very easy to give. For this, we pull back $P$ over the universal covering $\Sigma M$ of $M$, which is simply connected. This gives a principal bundle $\tilde P$, with structure group $G$ and with basis $\Sigma M$, the universal covering of $M$, where the connection $\theta$ can be pulled back into an unique connection $\tilde \theta$.

\begin{Definition} Given a category $\frak A$ of Lie groups. Let $\frak B$ be the set of objects in $\frak A$ such that,  for each group $H \in \frak B$ with Lie algebra $\frak h$, there is an injective morphism of Lie groups $\rho : H \rightarrow G$ such that the connection $\theta$ with curvature $\Omega$ is such that: 

for any smooth 1-parameter family $Hc_t$ of horizontal paths starting at $p$, for any smooth vector fields $X,Y$ in $TM$,  
\begin{equation*} s, t \in [0,1]^2 \rightarrow \Omega_{c_t(s)}(X,Y) \end{equation*} 
is a smooth $\frak h$-valued map,
\vskip 5pt \noindent
if the set $\frak B$ has a minimal element for the inclusion $H_0^{red}$, we call it the restricted holonomy group of class $\frak A$ of $\theta$. 

\end{Definition}

\begin{Theorem}\label{Holonomie}
In these three contexts:

\begin{item} 1) $G$ is CBH and $\frak A$ is the set of CBH Lie groups that are closed subgroups in $G$, \end{item}

\begin{item} 2) $G$ is of the first type and $\frak A$ is the set of Lie groups of the first type that are closed subgroups in $G$, \end{item}

\begin{item} 3) $G$ is of the second type and $\frak A$ is the set of Lie groups of the second type that are closed subgroups in $G$, \end{item} 

\vskip 5pt
\noindent
we can define the restricted holonomy group $H_0^{red}$ of class $\frak A$.Moreover, the Lie algebra of $H_0^{red}$ is the closed Lie algebra contained in $\frak g$ generated by the curvature elements. \end{Theorem}

\noindent
\textbf{Proof of the theorem \ref{Holonomie}:}
Since we work on the bundle $\tilde P$ that has a simply connected basis, we can restrict ourselves to the Lie groups of the set $\frak B$ deduced from $\frak A$ that are connected. 
 Then, using Theorem \ref{LieII} and the lemma following, we have that $\frak B$ has a minimal element, which is the connected Lie group with algebra spanned by the curvature elements. 

 The condition of closedness implies the enlargeability in these three cases. Then, since the Lie algebra of $H_0^{red}$ needs to contain all the curvature elements, the Lie algebra which is given is clearly the smaller one. \qed

\section{The holonomy group}

The holonomy group $H^{curv}$ induced by horizontal paths is not actually an infinite dimensional Lie group, principaly because its structure is induced by the structure of the non associative monoid of the loops based at $x$. This is why we need to define a Lie group $H$ such that there is morphism of topological group $ H^{curv} \rightarrow H$. $H$ needs to be the smaller as possible.

\begin{Definition}   
Let $\frak A$ be a category of Lie groups. We define the Holonomy group of class $\frak A$ as the minimal element of the set $\frak B'$ of objects $H$ with Lie algebra $\frak h$ of $\frak A$ such that
the connection $\theta$ with curvature $\Omega$ is such that

\noindent - for any smooth 1-parameter family $Hc_t$ of horizontal paths starting at $p$, for any smooth vector fields $X,Y$ in $TM$,  
\begin{equation*} s, t \in [0,1]^2 \rightarrow \Omega_{c_t(s)}(X,Y)  \end{equation*} 
is a smooth $\frak h$-valued map,

\noindent - $\frak h$ is stable under the action of $H^{curv}$.

\vskip 5pt \noindent
Then, if $\frak B'$ has a minimal element $H^{red}$, we call it Holonomy group of class $\frak A$ of the connection $\theta$.
 
\end{Definition}

With the conditions imposed to the groups of $\frak B '$, one can recognize the conditions of Theorem \ref{CourbureII} and of Lemma \ref{semidirect}, and hence see that these conditions ensure that $\theta$ reduces to the holonomy group if it is defined.

\begin{Theorem}\label{HolonomieII}
In these three contexts:

\begin{item} 1) $G$ is CBH and $\frak A$ is the set of CBH Lie groups that are closed subgroups in $G$, \end{item}

\begin{item} 2) $G$ is of the first type and $\frak A$ is the set of Lie groups of the first type that are closed subgroups in $G$, \end{item}

\begin{item} 3) $G$ is of the second type and $\frak A$ is the set of Lie groups of the second type that are closed subgroups in $G$, \end{item} 

\vskip 5pt \noindent
we can define the holonomy group $H^{red}$ of class $\frak A$. Moreover, the Lie algebra of $H$ is spanned by the curvature elements. \end{Theorem}

\noindent
\textbf{Proof of the Theorem \ref{HolonomieII}:}
According to  Theorem \ref{Holonomie} and Lemma \ref{semidirect}, we have only to show that the Lie algebra spanned by the curvature elements is stable by $H^{curv}$. We already know that it is stable by $H^{curv}_0$. Let $\{c_i\}$ be a family of loops generating $\pi_1(M)$, with holonomies $\{g_i\}$. Let $\Omega(X,Y)$ be a curvature element. Then, there exists a one parameter family of smooth loops $c_t$ based at $x$ such that $Hc_0(1)=Hc_0(0)$ and if $Hc_t(1)=Hc_t(0).g_t$, $\partial_{t=0}g_t = \Omega(X,Y)$. Hence, if $*$ is the mono\"id law of the loops based at $x$, we have that $(c_i * c_t ) * c_i^{-1}$ is a $0$-homotopic loop. Hence, its holonomy is in $ H^{curv}_0 \subset H_0^{red}$. Hence, taking the derivative with respect to $t$, we have that $Ad_{g_i}(\Omega(X,Y))$ is in the Lie algebra of $H_0^{red}$.  \qed

\section{Final remarks: the holonomy group as an invariant of the connection}
The holonomy groups $H^{red}$ considered here are closed subgroups of the structure group $G$. That is, they depend on $G$ as well as on $\theta$ and on  the basis manifold $M$. If we change the structure group $G$ for another Lie group with weaker structure, we also change the holonomy group. 

Conceptually, the holonomy group should depend only on the choice of the connection $\theta$ and on the basis manifold $M$. It is the case for $H^{curv}$, whose topology is the push-forward of the topology of the loops on $M$. One could expect to find a wider class of Lie subgroups  of $G$ than the class of closed Lie subgroups that enables one to define $H^{red}$ with a stronger topology than the one of $G$.  This leads to the following question:

What is the weaker property $(P)$ on Lie subalgebras of $\frak g$ such that 

- if a totally ordered family $B$ of subalgebras of $\frak g$  satisfy $(P)$, their intesection (equipped with the projective topology of the family $B$) also satisfies $(P)$,

- if $\frak h_1$ and $\frak h_2$ satisties $(P)$, then $\frak h_1 \cap \frak h_2$ satisfies $(P),$

- if $\frak h$ satisfies $(P)$, there is a Lie subgroup $H$ of $G$ with Lie algebra $\frak h$ ?

Under these assumptions, one could apply Zorn's lemma to find the holonomy group $H^{red}$ of the considered connection with respect to the class of Lie groups defined by $(P)$. The topology of this holonomy group should depend less on the choice on the structure group $G$.  

\vskip 10pt

More generally, a natural question arises on $H^{curv}$: when is it a Lie subgroup of $G$ ? This question is the same as: when is $H^{curv}_0$ a Lie subgroup of $G$ ? This question is related to the actual loss of easy criteria to decide whether a topological group is a (regular ?) infinite dimensional Lie group. If $H^{curv}$ is a regular infinite dimensional Lie group modelled on a complete locally convex topological space, then the theorem of reduction can be applied, and one can expect to find a category $\frak A$ of Lie groups for which $H^{curv}=H^{red}$. Otherwise, and especially if $H^{curv}$ is only an infinite dimensional Lie group (non regular) or only a topological group, the theorem of reduction cannot be applied and one can expect to have $H^{curv} \neq H^{red}$.
 

\end{document}